\documentclass[12pt]{article}
\usepackage{amsmath,amssymb,amsthm}
\usepackage[margin=1in]{geometry}
\usepackage{hyperref}
\usepackage{tikz}
\usepackage{graphicx}
\usepackage{url}
\usepackage[mathlines]{lineno}
\usepackage{multirow}
\usepackage{dsfont} % needed for double-struck 1
\usepackage{color}
\definecolor{red}{rgb}{1,0,0}

\newtheorem{thm}{Theorem}[section]
\newtheorem{cor}[thm]{Corollary}
\newtheorem{lem}[thm]{Lemma}

\newtheorem{obs}[thm]{Observation}

\theoremstyle{definition}
\newtheorem{defn}[thm]{Definition}

\theoremstyle{example}

\def\mtx#1{\begin{bmatrix} #1 \end{bmatrix}}

\DeclareMathOperator{\chr}{char}

\DeclareMathOperator{\pr}{pr}
\DeclareMathOperator{\epr}{epr}

\DeclareMathOperator{\rank}{rank}

\newcommand{\Fnn}{F^{n\times n}}

\newcommand{\x}{\times}
\newcommand{\bit}{\begin{itemize}}
\newcommand{\eit}{\end{itemize}}
\newcommand{\ben}{\begin{enumerate}}
\newcommand{\een}{\end{enumerate}}
\newcommand{\beq}{\begin{equation}}
\newcommand{\eeq}{\end{equation}}
\newcommand{\bea}{\begin{eqnarray*}}
\newcommand{\eea}{\end{eqnarray*}}
\newcommand{\bpf}{\begin{proof}}
\newcommand{\epf}{\end{proof}\ms}
\newcommand{\bmt}{\begin{bmatrix}}
\newcommand{\emt}{\end{bmatrix}}
\newcommand{\ms}{\medskip}

\newcommand{\ba}{\begin{array}}
\newcommand{\ea}{\end{array}}

\usepackage{enumitem}
\DeclareMathOperator{\qpr}{qpr}
\DeclareMathOperator{\SPAN}{span}

\begin{document}
\title{The quasi principal rank characteristic sequence}
\author{
Shaun M. Fallat\thanks{Department of Mathematics and Statistics,
University of Regina,  Regina, Saskatchewan, S4S0A2, Canada. Research supported in part by an NSERC Discovery Research Grant.
V8W 2Y2, Canada (shaun.fallat@uregina.ca).}
\and
Xavier Mart\'inez-Rivera\thanks{Department of Mathematics and Statistics,
Auburn University, Auburn, AL 36849, USA (martinez.rivera.xavier@gmail.com).}
}

%\linenumbers

\maketitle

\begin{abstract}
A minor of a matrix is quasi-principal if it is a principal or an almost-principal minor.
The quasi principal rank characteristic sequence (qpr-sequence) of
an  $n\times n$ symmetric matrix is introduced, which is defined as
$q_1 q_2 \cdots q_n$, where $q_k$ is  {\tt A}, {\tt S}, or {\tt N},
according as all, some but not all, or none of its quasi-principal minors of
order $k$ are nonzero.
This sequence extends the principal rank characteristic sequences in
the literature, which only depend on the principal minors of the matrix.
A necessary condition for the attainability of a qpr-sequence is established.
Using probabilistic techniques,
a complete characterization of the qpr-sequences that are attainable by
symmetric matrices over fields of characteristic $0$ is given.
\end{abstract}


\medskip


%%%%%%%%%%%%%%%%%%%%%%%%%%%%%%%%%%%%%%%
%%%%%%%%%%%%%%%%%%%%%%%%%%%%%%%%%%%%%%%
%%%%%%%%%%%%%%%%%%%%%%%%%%%%%%%%%%%%%%%
%%%%%%%%%%%%%%%%%%%%%%%%%%%%%%%%%%%%%%%
%%%%%%%%%%%%%%%%%%%%%%%%%%%%%%%%%%%%%%%
%%%%%%%%%%%%%%%%%%%%%%%%%%%%%%%%%%%%%%%
%%%%%%%%%%%%%%%%%%%%%%%%%%%%%%%%%%%%%%%
%%%%%%%%%%%%%%%%%%%%%%%%%%%%%%%%%%%%%%%
%\pagebreak
\section{Introduction}\label{sintro}
$\null$
\indent
Motivated by the famous principal minor assignment problem,
which is stated in \cite{HS}, Brualdi et al.\ \cite{P} introduced the principal rank characteristic sequence,
which is defined as follows:
Given an $n \times n$ symmetric matrix $B$ over a
field $F$,
\textit{the principal rank characteristic sequence}
(abbreviated pr-sequence) of $B$ is
$\pr(B) = r_0]r_1 \cdots r_n$, where, for $k \geq 1$,
   \begin{equation*}
      r_k =
         \begin{cases}
             1 &\text{if $B$ has a nonzero principal minor of order $k$, and}\\
             0 &\text{otherwise,}
         \end{cases}
   \end{equation*}
while $r_0 = 1$ if and only if $B$ has a
$0$ diagonal entry.
(The \textit{order} of a minor is $k$ if it is
the determinant of a $k \times k$ submatrix.)

The pr-sequence was later ``enhanced'' by Butler et al.\ \cite{EPR} with the introduction of another sequence:
Given an $n \times n$ symmetric matrix $B$ over a
field $F$,
the \textit{enhanced principal rank characteristic sequence} (abbreviated epr-sequence) of $B$ is
$\epr(B) = \ell_1\ell_2 \cdots \ell_n$, where
   \begin{equation*}
      \ell_k =
         \begin{cases}
             \tt{A} &\text{if all the principal minors of order $k$ are nonzero;}\\
             \tt{S} &\text{if some but not all of the principal minors of order $k$ are nonzero;}\\
             \tt{N} &\text{if none of the principal minors of order $k$ are nonzero, i.e., all are zero.}
         \end{cases}
   \end{equation*}
Pr- and epr-sequences have received considerable attention since their introduction
(see \cite{{P}, {BIRS13}, {EPR}, {skew}, {XMR-Classif}, {EPR-Hermitian}, {XMR-Char 2}}, for example), and, recently, another sequence, one that refines the epr-sequence,
called the signed enhanced principal rank characteristic sequence,
was introduced by Mart\'inez-Rivera in \cite{XMR-sepr}.
A commonality among these sequences is the fact that the only minors about which they provide information are the principal minors of the matrix.
Thus, a natural problem is to determine what epr-sequences are attainable if the definition of the epr-sequence is extended by expanding the set of minors it considers; that is, the aim of this paper.
The question is, which other minors of the matrix should we consider?
Our choice, which we justify below, is the set of almost-principal minors.
This leads us to introduce a new sequence,
which will be the focus of this paper.
Before formally defining this sequence, some terminology is necessary:
For a given $n \times n$ matrix $B$ and
$\alpha, \beta \subseteq \{1,2, \dots, n\}$,
$B[\alpha, \beta]$ denotes the submatrix lying in
rows indexed by $\alpha$ and columns indexed by $\beta$;
$B[\alpha, \beta]$ is called a {\em quasi-principal} submatrix if
$|\alpha| = |\beta|$ and
$|\alpha| - 1  \leq  |\alpha \cap \beta| \leq |\alpha|$;
we call the minor $\det B[\alpha, \beta]$ {\em quasi-principal} if
$B[\alpha, \beta]$ is a  quasi-principal submatrix, and we say that
its \emph{order} is $|\alpha|$.
Note that a minor is quasi-principal if
it is a principal or an almost-principal minor.
In what follows, $F$ is used to denote a field.

\begin{defn}{\rm
The {\em quasi principal rank characteristic sequence} of a symmetric matrix $B\in\Fnn$ is
the sequence (qpr-sequence)  $\qpr(B)=q_1q_2\cdots q_n$, where
\[q_k=\left\{\begin{array}{ll}
{\tt A} &  \mbox{ if  all of the quasi-principal minors of order $k$ are nonzero;}\\
{\tt S} &   \mbox{ if some but not all of the quasi-principal minors of order $k$ are nonzero;}\\
{\tt N} &  \mbox{ if none of the quasi-principal minors of order $k$ are nonzero, i.e., all are zero.}\end{array}\right.\]
}\end{defn}

There are several reasons why the almost-principal minors are a natural choice of minors for extending the epr-sequence with the qpr-sequence.
For example, almost-principal minors (and principal minors) are of interest in algebraic geometry and statistics:
In \cite{Stu09}, Sturmfels  presented mathematical problems whose solutions would, according to him, ``likely become important contributions to the emerging interactions between algebraic geometry and computational statistics.''
One of the problems he presented is to study the geometry of conditional independence models for multivariate Gaussian random variables, where
almost-principal minors play a crucial role, due to their
intimate connection with conditional independence in statistics \cite{Stu09};
the program presented in \cite{Stu09} was recently continued by
Boege, D'Al\`{i}, Kahle and Sturmfels
in \cite{Stu17}.
We note that almost-principal minors are referred to as
\textit{partial covariance}
(or, if renormalized, \textit{partial correlations}) in
the statistics literature \cite{Stu09}.

Almost-principal minors also play a role in theoretical physics:
Motivated by the fact that the Yang-Baxter equation for certain models  related to free fermions in statistical mechanics arise from algebraic identities among the minors of the corresponding matrix kernel,
Kenyon and Pemantle \cite{KP14} studied algebraic relations between the principal minors of a matrix.
They wanted to understand in what sense these identities are the only algebraic identities among the correlations.
By considering almost-principal minors, they showed that a complete description in terms of the translates of a single relation (the hexahedron relation) can be obtained.
Related work was subsequently done by Sturmfels et al.\ in \cite{Stu16}, where a conjecture of Kenyon and Pemantle \cite{KP14} was proven.

Almost-principal minors are also relevant in
the study of $M$-matrices, that is,
real matrices having
non-positive off-diagonal entries and
eigenvalues with positive real part.
%A real matrix is a $Z$-matrix if its off-diagonal entries are nonpositive.
%A real matrix is an $M$-matrix if it is a $Z$-matrix having eigenvalues with positive real part.
A matrix that occurs as the inverse of an $M$-matrix is called an \textit{inverse $M$-matrix}.
A question of interest in some applications of $M$-matrices is the following: When is a nonnegative matrix an inverse $M$-matrix?
The reader is referred to
\cite{InvM1} and the surveys in \cite{{InvM2},{InvM3}} for
more information.
%\begin{quest}\label{Quest: Z-matrix}
%For nonsingular (entrywise) nonnegative matrices $B$,
%when is $B^{-1}$ a $Z$-matrix?
%\end{quest}
The relevance of almost-principal minors (and principal minors) in this context becomes evident when we note that,
for a given nonsingular matrix $B$,
each entry of $B^{-1}$ can be expressed in
terms of ratios of quasi-principal minors of $B$.

The remainder of this section is devoted to
list some known results that will be used frequently in what follows.
In Section \ref{s: qpr-sequence},
we state some basic observations about qpr-sequences, and
conclude with a necessary condition for the attainability of a qpr-sequence.
In Section \ref{s: Main results},
we confine our attention to the attainable qpr-sequences of
symmetric matrices over fields of characteristic 0,
which we completely characterize in our main result,
Theorem \ref{thm: main result}.

The qpr-sequence of $n \times n$ matrix is said to have \textit{length} $n$.
Given a vector $x$ of length $n$,
$x[\alpha]$ denotes the subvector of $x$ with
entries indexed by $\alpha \subseteq \{1,2, \dots, n\}$.
The column and row spaces of a matrix $B$ are denoted by
$C(B)$ and $R(B)$, respectively.

%%%%%%%%%%%%%%%%%%%%%%%%%%%%%%%%%%%%%%%
%%%%%%%%%%%%%%%%%%%%%%%%%%%%%%%%%%%%%%%
%%%%%%%%%%%%%%%%%%%%%%%%%%%%%%%%%%%%%%%
%%%%%%%%%%%%%%%%%%%%%%%%%%%%%%%%%%%%%%%
%%%%%%%%%%%%%%%%%%%%%%%%%%%%%%%%%%%%%%%
%%%%%%%%%%%%%%%%%%%%%%%%%%%%%%%%%%%%%%%
%%%%%%%%%%%%%%%%%%%%%%%%%%%%%%%%%%%%%%%
%%%%%%%%%%%%%%%%%%%%%%%%%%%%%%%%%%%%%%%
\subsection{Previous works}
$\null$
\indent
In this section, we list known results that
will be used frequently in subsequent sections.
\begin{thm}
{\rm \cite[Theorem 2.3]{EPR}}
{\rm ({\tt NN} Theorem.)}
\label{thm:NN}
Let $B \in \Fnn$ be symmetric.
Suppose that
$\epr(B) = \ell_1 \ell_2 \cdots \ell_n$ and that $\ell_k = \ell_{k+1} = {\tt N}$ for some $k$.
Then $\ell_j = {\tt N}$ for all $j \geq k$.
\end{thm}

The following well-known fact is stated in \cite{BIRS13}.

\begin{thm}
\label{thm: rank of a symm mtx}
{\rm \cite[Theorem 1.1]{BIRS13}}
Let $B \in \Fnn$ be symmetric.
Then $\rank B = \max \{ |\gamma| : \det (B[\gamma]) \neq 0 \}$
(where the maximum over the empty set is defined to be 0).
\end{thm}

Another well-known fact is stated next.

\begin{lem}\label{rank when appending}
Let
\[
B = \mtx{B[\alpha] & x \\
        y^T & b} \in \Fnn
\]
and  $|\alpha| = n-1$.
Then $\rank B[\alpha] \leq \rank B \leq \rank B[\alpha] +2$.
Furthermore, the following hold:
\ben
\item[(i)] If $x \notin C(B[\alpha])$ and $y^T \notin R(B[\alpha])$, then
$rank B = \rank B[\alpha] +2$.

\item[(ii)]
If $x \in C(B[\alpha])$ and $y^T \notin R(B[\alpha])$, then
$rank B = \rank B[\alpha] +1$.

\item[(iii)]
If $x \notin C(B[\alpha])$ and $y^T \in R(B[\alpha])$, then
$rank B = \rank B[\alpha] +1$.

\item[(iv)] If $x \in C(B[\alpha])$ and $y^T \in R(B[\alpha])$, then
$\rank B[\alpha] \leq \rank B \leq \rank B[\alpha] +1$,
and
$\rank B = \rank B[\alpha] + 1$
if and only if
$b \neq v^T B[\alpha] u$, where
$u$ and $v$ are such that
$x=B[\alpha]u$ and $y=B[\alpha]v$.
\een
\end{lem}

In the interest of clarity, we state a special case of
Lemma \ref{rank when appending} that will be used repeatedly.

\begin{lem}\label{rank when appending same row and column}
Let
\[
B = \mtx{B[\alpha] & y \\
        y^T & b} \in \Fnn
\]
and $|\alpha| = n-1$.
Then $\rank B[\alpha] \leq \rank B \leq \rank B[\alpha] +2$.
Furthermore, the following hold:
\ben
\item[(i)] $rank B = \rank B[\alpha] +2$ if and only if
$y \notin C(B[\alpha])$.
\item[(ii)] $\rank B = \rank B[\alpha] + 1$ if and only if
$y \in C(B[\alpha])$ and  $b \neq x^TB[\alpha]x$,
where $x$ is such that $y=B[\alpha]x$.
\een
\end{lem}

The next lemma follows from Lemma \ref{rank when appending}.

\begin{lem}\label{rank when deleting}
Let $B \in \Fnn$ be symmetric and nonsingular.
Then the rank of any $(n-1)\times (n-1)$ submatrix of $B$ is at least $n-2$.
\end{lem}

%%%%%%%%%%%%%%%%%%%%%%%%%%%%%%%%%%%%%%%
%%%%%%%%%%%%%%%%%%%%%%%%%%%%%%%%%%%%%%%
%%%%%%%%%%%%%%%%%%%%%%%%%%%%%%%%%%%%%%%
%%%%%%%%%%%%%%%%%%%%%%%%%%%%%%%%%%%%%%%
%%%%%%%%%%%%%%%%%%%%%%%%%%%%%%%%%%%%%%%
%%%%%%%%%%%%%%%%%%%%%%%%%%%%%%%%%%%%%%%
%%%%%%%%%%%%%%%%%%%%%%%%%%%%%%%%%%%%%%%
%%%%%%%%%%%%%%%%%%%%%%%%%%%%%%%%%%%%%%%
\section{The quasi principal rank characteristic sequence}
\label{s: qpr-sequence}
$\null$
\indent
This section is devoted to stating some basic observations about qpr-sequences, and to
establish a necessary condition for the attainability of a qpr-sequence.

\begin{obs}
\label{obs: rank and A or N at end}
Let $B \in \Fnn$ be symmetric.
Then $\qpr(B)$ must end in ${\tt N}$ or ${\tt A}$, and
$\rank B$ is equal to the index of the last {\tt A} or {\tt S} in qpr$(B)$.
\end{obs}

\begin{obs}
\label{obs: qpr with A and N}
Let $B \in \Fnn$ be symmetric,
let $\qpr(B)=q_1q_2 \cdots q_n$, and
let $\epr(B)= \ell_1\ell_2 \cdots \ell_n$.
If $q_k \in \{\tt A, N\}$ for some $k$, then $\ell_k = q_k$.
\end{obs}

\begin{obs}
\label{obs: N at start}
Let $B \in \Fnn$ be symmetric.
Suppose that $\qpr(B) = q_1 q_2 \cdots q_n$ and that
$q_1 = \tt N$.
Then $B$ is the zero matrix and $q_j = \tt N$ for all $j \geq 1$.
\end{obs}

We now generalize Observation \ref{obs: N at start}
in Theorem \ref{thm:N implies all N} below;
but before that, we need the tools provided by the next two results.
Given a matrix $B$ with a nonsingular principal
submatrix $B[\gamma]$, we denote by $B/B[\gamma]$ the
Schur complement of $B[\gamma]$ in $B$ (see \cite{Schur}).
The following result is a well-known fact.

\begin{thm}
\label{schur}
{\rm (Schur Complement Theorem.)}
Let $B \in \Fnn$ be symmetric with
$\rank B=r$.
Let $B[\gamma]$ be a nonsingular principal
submatrix of $B$ with $|\gamma| = k \leq r$,
and let $C = B/B[\gamma]$.
Then the following results hold: %\vspace{-3mm}
\begin{enumerate}
\item [$(i)$]\label{p1SC}
$C$ is an $(n-k)\times (n-k)$ symmetric matrix. %\vspace{-3mm}
\item [$(ii)$]\label{p2SC}
Assuming that the indexing of $C$ is inherited from $B$,
any minor of $C$ is given by %\vspace{-3mm}
\[ \det C[\alpha, \beta] = \det B[\alpha \cup \gamma, \beta \cup \gamma]/ \det B[\gamma].\]
\item [$(iii)$]\label{p3SC} $\rank C = r-k$.%\vspace{-3mm}
%\item [$(iv)$]\label{LHp4SC} Any nonsingular principal submatrix of $B$ of
%order at most $r$ is contained in a nonsingular principal submatrix of order $r$.
\end{enumerate}
\end{thm}

A corollary now follows immediately from the Schur Complement Theorem:

\begin{cor}\label{schurAN}
{\rm (Schur Complement Corollary.)}
Let $B \in \Fnn$ be symmetric,
let $\qpr(B)=q_1 q_2 \cdots q_n$, and
let $B[\gamma]$ be a nonsingular
principal submatrix of $B$,
with $|\gamma| = k \leq \rank B$.
Let $C = B/B[\gamma]$
and $\qpr(C)=q'_{1} q'_2 \cdots q'_{n-k}$.
Then, for $j=1, \dots, n-k$,
$q'_j=q_{j+k}$  if
$q_{j+k} \in \{{\tt A,N}\}$.
\end{cor}

\begin{thm}
\label{thm:N implies all N}
{\rm ({\tt N} Theorem.)}
Let $B \in \Fnn$ be symmetric.
Suppose that $\qpr(B) = q_1 q_2 \cdots q_n$ and that
$q_k = \tt N$ for some $k$.
Then $q_j = \tt N$ for all $j \geq k$.
\end{thm}

\bpf
If $k=1$, the desired conclusion follows from Observation \ref{obs: N at start}.
Now assume that $k \geq 2$ and
suppose to the contrary that $q_j \neq \tt N$ for some $j \geq k$.
Hence, $\rank B \geq j$.
Let $\epr(B) = \ell_1 \ell_2 \cdots \ell_n$.
Since $q_k = \tt N$, $\ell_k = \tt N$.
Note that $\ell_{k-1} \neq \tt N$, as otherwise
the $\tt NN$ Theorem would imply that $\ell_j = \tt N$ for $j > k-1$,
and therefore that $\rank B < k-1 < j$
(see Theorem \ref{thm: rank of a symm mtx}).
Let $B[\gamma]$ be a $(k-1) \times (k-1)$ nonsingular principal submatrix.
By the Schur Complement Theorem, the rank of $B/B[\gamma]$ is
$\rank B - (k-1) \geq j -(k-1) > 0$.
Moreover, by the Schur Complement Corollary,
the qpr-sequence of $B/B[\gamma]$ begins with $\tt N$,
which implies that $B/B[\gamma]$ is the zero matrix
(see  Observation \ref{obs: N at start}), whose rank is 0,
leading to a contradiction.
\epf

\begin{cor}
\label{cor: necessary condition}
If the qpr-sequence $q_1q_2 \cdots q_n$ is attainable by a
symmetric matrix over a field $F$, then the following statements hold:
\ben
\item[(i)] $q_n \neq \tt S$.
\item[(ii)] Neither $\tt NA$ nor $\tt NS$ is a subsequence of $q_1q_2 \cdots q_n$.
\een
\end{cor}

A natural question now arises:
Is the necessary condition in Corollary \ref{cor: necessary condition} also sufficient?
For fields of characteristic zero,
this question is answered affirmatively in the next section.

%%%%%%%%%%%%%%%%%%%%%%%%%%%%%%%%%%%%%%%
%%%%%%%%%%%%%%%%%%%%%%%%%%%%%%%%%%%%%%%
%%%%%%%%%%%%%%%%%%%%%%%%%%%%%%%%%%%%%%%
%%%%%%%%%%%%%%%%%%%%%%%%%%%%%%%%%%%%%%%
%%%%%%%%%%%%%%%%%%%%%%%%%%%%%%%%%%%%%%%
%%%%%%%%%%%%%%%%%%%%%%%%%%%%%%%%%%%%%%%
%%%%%%%%%%%%%%%%%%%%%%%%%%%%%%%%%%%%%%%
%%%%%%%%%%%%%%%%%%%%%%%%%%%%%%%%%%%%%%%
\section{Main results and probabilistic techniques}
\label{s: Main results}
$\null$
\indent
In this section, using certain probabilistic techniques,
a complete characterization of the qpr-sequences that
are attainable by symmetric matrices over fields of characteristic 0 is established (see Theorem \ref{thm: main result});
this is done by demonstrating that the necessary condition established in
Corollary \ref{cor: necessary condition} for the attainability of a
qpr-sequence is  also sufficient, when the field is of characteristic $0$.

\begin{lem}
\label{copyN}
Let $F$ be a field with $\chr F=0$, and let  $B\in\Fnn$ be symmetric.
Assume  that  $\qpr(B)=q_1q_2 \cdots q_n$ and  $r=\rank B$.
Let
\[
B' = \mtx{B & Bx \\
                x^T B & x^TBx},
\]
where $x$ is chosen randomly, and
let $\qpr(B') = q_1'q_2' \cdots q_n'q_{n+1}'$.
Then
\ben
\item[(i)] $q_i'={\tt N}$ for $i=r+1,\dots,n+1$.
\item[(ii)] For $1 \le i \le r$,
$q_i'={\tt S}$ if $q_i={\tt S}$ or %$q_i={\tt N}$.
\item[(iii)] For $1 \le i \le r$,
with high probability,  $q_i'={\tt A}$ if $q_i={\tt A}$.
\een
\end{lem}
\bpf
First, note that, since $\rank(B) =r$,
the $\tt N$ Theorem implies that $q_i \neq \tt N$ for $1 \leq i \leq r$
(see Observation \ref{obs: rank and A or N at end}).
Hence, $q_j \in \{\tt A,S \}$ for $j=1,2, \dots, r$.

To see why (i) holds, observe that $\rank(B') = \rank(B) = r$,
which implies that $q'_i = \tt N$ for $i=r+1,\dots,n+1$.

To establish (ii), note that
the desired conclusion is immediate when $q_i = \tt S$.
%when $q_i = \tt N$, the conclusion follows vacuously,
%since $q_j \in \{\tt A,S \}$ for $j=1,2, \dots, r$.

We now show that (iii) holds.
Suppose that $q_i = \tt A$.
Let $C$ be an $i \times i$ quasi-principal submatrix of $B'$.
If $C$ is also a submatrix of $B$, then
$C$ is nonsingular by the assumption that $q_i = \tt A$.
Hence, suppose $C$ is not a submatrix of $B$,
meaning that $C$ involves the last row or last column of $B'$.
We proceed by considering two cases.

\noindent
{\bf Case 1}: $C$ is a principal submatrix of $B'$.

\noindent
Then
\[
C = \mtx{C[\alpha] & y \\
                y^T  & x^TBx}
\]
for some $\alpha \subseteq \{1,2, \dots, n\}$ with $|\alpha| = i-1 <r$,
where $y = Bx[\alpha]$.
Since $C[\alpha]$ is an $(i-1) \times (i-1)$ submatrix of
an $i \times i$ principal submatrix of $B$,
and because $q_i = \tt A$,
Lemma \ref{rank when deleting} implies that $\rank C[\alpha] \geq i-2$.
We now examine two cases based on the value of $\rank C[\alpha]$.

\noindent
{\bf Subcase A}: $\rank C[\alpha] = i-2$.

\noindent
Then, with high probability, $x$ was chosen so that
$y$ is not in the range of $C[\alpha]$,
implying that $\rank C = \rank C[\alpha] + 2 = i$
(see Lemma \ref{rank when appending same row and column}).
Since $C$ is an $i \times i$ matrix, $C$ is nonsingular, as desired.

\noindent
{\bf Subcase B}: $\rank C[\alpha] = i-1$.

\noindent
Then $C[\alpha]$ has full rank, and therefore
$y$ is in the column space of $C[\alpha]$.
With high probability,
$x$ was chosen so that $(x[\alpha])^TC[\alpha](x[\alpha]) \neq x^TBx$;
hence,  $\rank C = \rank C[\alpha] + 1 = i$
(see Lemma \ref{rank when appending same row and column}),
implying that $C$ is nonsingular, as desired.

\medskip
\medskip
\noindent
{\bf Case 2}: $C$ is an almost-principal submatrix of $B'$ involving
the last row or last column of $B'$, or both.
Since $B'$ is symmetric, it suffices to examine
the following two subcases of Case 2.

\noindent
{\bf Subcase A}: $C$ involves the last column of $B'$ but not the last row of $B'$.

\noindent
Suppose that the index of the last row of $C$ is $k$.
It is easy to see that
\[
C = \mtx{B[\alpha] & y \\
                z^T  & B[\{k\},\{n+1\}]}
\]
for some $y$ and $z$,
and some $\alpha \subseteq \{1,2, \dots, n\}$ with $|\alpha| = i-1 <r$.
Since the first $|\alpha|=i-1$ columns of $C$ are columns of
the matrix $B[\alpha \cup \{k\}]$, which is nonsingular (because $q_i = \tt A$), the matrix
\[
\mtx{B[\alpha] \\ z^T}
\]
has full rank.
With high probability, $x$ was chosen so that the vector
\[
\mtx{y \\ B[\{k\},\{n+1\}]}
\]
is not in the span of the first $i-1$ columns of $C$;
hence, $C$ is nonsingular, as desired.

\noindent
{\bf Subcase B}: $C$ involves both the last row and the last column of $B'$.

\noindent
Then
\[
C = \mtx{B[\alpha, \beta] & y \\
                z^T  & x^TBx}
\]
for some $\alpha, \beta \subseteq \{1,2, \dots, n\}$ with
$|\alpha| = |\beta| = i-1 <r$ and $|\alpha \cap \beta| = |\alpha|-1$,
where $y = Bx[\alpha]$ and $z = Bx[\beta]$.
Since $B[\alpha, \beta]$ is an $(i-1) \times (i-1)$ submatrix of
an $i \times i$ quasi-principal submatrix of $B$,
and because $q_i = \tt A$,
Lemma \ref{rank when deleting} implies that
$\rank B[\alpha, \beta] \geq i-2$.

\noindent
{\bf Sub-subcase a}: $\rank B[\alpha, \beta] = i-2$.

\noindent
Then with high probability $x$ was chosen so that
$y$ is not in the column space of $B[\alpha, \beta]$ and
$z^T$ is not in the row space of $B[\alpha, \beta]$,
implying that
$\rank C = \rank B[\alpha, \beta]+2 = i$
(see Lemma \ref{rank when appending}), and therefore that
$C$ is nonsingular, as desired.

\noindent
{\bf Sub-subcase b}: $\rank B[\alpha, \beta] = i-1$.

\noindent
In this case, $\rank B[\alpha, \beta]$ has full rank,
meaning that
$y$ is in the column space of $B[\alpha, \beta]$, and that
$z^T$ is in the row space of $B[\alpha, \beta]$.
With high probability, $x$ was chosen so that
the rank of $C$ is larger than that of $B[\alpha, \beta]$, implying that
$\rank C = \rank B[\alpha, \beta]+1 = i$
(see Lemma \ref{rank when appending}), and therefore that
$C$ is nonsingular, as desired.
\epf

%From this point on, given the above established claim, we assume that $q_i \in \{ \tt A, \tt S\}$.

\begin{lem}
\label{copyA}
Let $F$ be a field with $\chr F=0$, and let  $B\in\Fnn$ be symmetric.
Assume that  $\qpr(B)=q_1q_2\cdots q_{n-1}{\tt A}$ with $q_k\in\{{\tt A, S}\}$ for
$k=1,\dots,n-1$.
Then there exists a symmetric matrix $B'\in F^{(n+1)\x(n+1)}$ such that
$\qpr(B')=q_1 \cdots q_{n-1}{\tt AA}$.
 \end{lem}

\bpf
Let
\[
C = \mtx{B & Bx \\
                x^T B & x^TBx},
\]
where $x$ is chosen randomly.
By Lemma \ref{copyN},
with high probability, $\qpr(C) = q_1q_2 \cdots q_{n-1} \tt AN$.
Let
\[
B' = \mtx{B & Bx \\
                x^T B & x^TBx + t},
\]
and suppose that $\qpr(B') = q'_1q'_2 \cdots q'_nq'_{n+1}$.
We now show that, with high probability,
$\qpr(B')=q_1 \cdots q_{n-1}{\tt AA}$.
Note that $\det B' = t \det B + \det C = t \det B$.
Then, as $\det B \neq 0$, with high probability, $t$ was chosen so that
$\det B' \neq 0$, implying that $q'_{n+1} = \tt A$.

Observe that if $q_j = \tt S$ for some $j$,
meaning that
$B$ contains both a zero and nonzero
quasi-principal minor of order $j$, then
$B'$ also does; hence,
if  $q_k = \tt S$ for some $k$ with $1 \leq k \leq n-1$, then
$q'_k = \tt S$.

By hypothesis, $q_n := \tt A$.
Now note that to conclude the proof, it suffices to show that
the following statement holds with high probability:
If $q_j = \tt A$ for some $j$ with $1 \leq j \leq n$,
then $q'_j = \tt A$.
Thus, suppose that $q_j = \tt A$ for some $j$ with $1 \leq j \leq n$.
Let $B'[\alpha, \beta]$ be a $j \times j$ quasi-principal submatrix of $B'$,
meaning that $|\alpha| = |\beta| = j$ and
$|\alpha| - 1  \leq  |\alpha \cap \beta| \leq |\alpha|$.
Since $q_j = \tt A$,
$\det C[\alpha, \beta] \neq 0$.
If $n+1 \notin \alpha \cup \beta$, then
$\det B'[\alpha, \beta]  = \det C[\alpha, \beta] \neq 0$, as desired.
Now suppose that $n+1 \in \alpha \cup \beta$.
Since $B'$ is symmetric, it suffices to examine the following two cases.

\noindent
{\bf Case 1}: $n+1 \in \alpha \cap \beta$.

\noindent
Then
$\det B'[\alpha, \beta] =
t \det B'[\alpha \setminus \{n+1\}, \beta  \setminus \{n+1\}] +
 \det C[\alpha, \beta]$.
Since $\det C[\alpha, \beta] \neq 0$, with high probability,
$t$ was chosen so that $\det B'[\alpha, \beta] \neq 0$.

\noindent
{\bf Case 2}: $n+1 \notin \alpha$ and $n+1 \in \beta$.

\noindent
Suppose that the index of the last row of $B'[\alpha, \beta]$ is $k$,
and that $l \in \alpha $ but $l \notin \beta$.
It is easy to see that
\[
B'[\alpha, \beta] = \mtx{B[\alpha \setminus \{k\}, \beta \setminus \{n+1\}] & y \\
                                            z^T  & B[\{k\},\{n+1\}]}
\]
for some vectors $y,z \in F^{j-1}$.

Since the first $|\alpha|-1=j-1$ columns of $B'[\alpha, \beta] $ are columns of the quasi-principal submatrix
$B[\alpha, \beta \cup \{l\} \setminus \{n+1\}]$,
which is nonsingular (because $q_j = \tt A$),
the matrix
\[
\mtx{B[\alpha \setminus \{k\}, \beta \setminus \{n+1\} ] \\ z^T}
\]
has full rank.
With high probability, $x$ was chosen so that the vector
\[
\mtx{y \\ B[\{k\},\{n+1\}]}
\]
is not in the span of the first $j-1$ columns of $B'[\alpha, \beta] $;
hence, $B'[\alpha, \beta] $ is nonsingular, as desired.
\epf

\begin{lem}\label{AtoSN}
Let $F$ be a field with $\chr F=0$, and let  $B\in\Fnn$ be symmetric.
Assume  that $\qpr(B)=q_1 q_2 \cdots q_{n-1}{\tt A}$ with
$\ell_k\in\{{\tt A, S}\}$  for $k=1,\dots,n$.
Then there exists a symmetric matrix $B'\in F^{(n+1)\x(n+1)}$ such that
$\qpr(B')=q_1 \cdots q_{n-1}{\tt SN}$.
 \end{lem}
\bpf
Suppose $B = [b_1, b_2, \dots, b_n]$, where $b_i$ is the $i$th column of $B$, and
let $z \in \SPAN\{b_2, b_3, \dots, b_n\}$ be chosen randomly,
and assume that $z=Bx$.
Let
\[
B' = \mtx{B & z \\
                z^T  & p},
\]
where $p = x^TBx$.
Let $\qpr(B') = q'_1q'_2 \cdots q'_n q'_{n+1}$.
We now show that
$\qpr(B') = q_1q_2 \cdots q_{n-1} \tt SN$.
By Lemma \ref{rank when appending same row and column},
$\rank B' = \rank B$, implying that $B'$ is singular;
hence, $q'_{n+1} = \tt N$.
Observe that
\[\det B'[\{1, \dots, n\}, \{2,  \dots, n+1\}] = 0 \mbox{ and }
\det B'[\{1, \dots, n\}] = \det B \neq 0,
\]
implying that
$q'_{n} = \tt S$.
Let $i$ be an integer with $1 \leq i \leq n-1$.
Since $B$ is a submatrix of $B'$, $q'_i = q_i$ if $q_i = \tt S$.
Now, suppose $q_i = \tt A$.
Let $C$ be an $i \times i$ quasi-principal submatrix of $B'$.
%%%%%%%%%%%%%%%%%%%%%%%%%%%%
%%%%%%%%%%%%%%%%%%%%%%%%%%%%
%%%%%%%%%%%%%%%%%%%%%%%%%%%%
%The proof of Cases 1 and 2 is done exactly as in Prop. 3.1.
%%%%%%%%%%%%%%%%%%%%%%%%%%%%
If $C$ is also a submatrix of $B$, then
$C$ is nonsingular by the assumption that $q_i = \tt A$.
Hence, suppose $C$ is not a submatrix of $B$,
meaning that $C$ involves the last row or last column of $B'$.
We proceed by considering two cases.

\noindent
{\bf Case 1}: $C$ is a principal submatrix of $B'$.

\noindent
Then
\[
C = \mtx{C[\alpha] & y \\
                y^T  & x^TBx}
\]
for some $\alpha \subseteq \{1,2, \dots, n\}$ with $|\alpha| = i-1 <r$,
where $y = Bx[\alpha]$.
Since $C[\alpha]$ is an $(i-1) \times (i-1)$ submatrix of
an $i \times i$ principal submatrix of $B$,
and because $q_i = \tt A$,
Lemma \ref{rank when deleting} implies that $\rank C[\alpha] \geq i-2$.
We now examine two cases based on the value of $\rank C[\alpha]$.

\noindent
{\bf Subcase A}: $\rank C[\alpha] = i-2$.

\noindent
Then, with high probability, $z$ was chosen so that
$y$ is not in the range of $C[\alpha]$,
implying that $\rank C = \rank C[\alpha] + 2 = i$
(see Lemma \ref{rank when appending same row and column}).
Since $C$ is an $i \times i$ matrix, $C$ is nonsingular, as desired.

\noindent
{\bf Subcase B}: $\rank C[\alpha] = i-1$.

\noindent
Then $C[\alpha]$ has full rank, and therefore
$y$ is in the column space of $C[\alpha]$.
With high probability,
$z$ was chosen so that $(x[\alpha])^TC[\alpha](x[\alpha]) \neq x^TBx$;
hence,  $\rank C = \rank C[\alpha] + 1 = i$
(see Lemma \ref{rank when appending same row and column}),
implying that $C$ is nonsingular, as desired.

\medskip
\medskip
\noindent
{\bf Case 2}: $C$ is an almost-principal submatrix of $B'$ involving
the last row or last column of $B'$, or both.

\noindent
Since $B'$ is symmetric, it suffices to examine
the following two subcases of Case 2.

\noindent
{\bf Subcase A}: $C$ involves the last column of $B'$ but not the last row of $B'$.

\noindent
Suppose that the index of the last row of $C$ is $k$.
It is easy to see that
\[
C = \mtx{B[\alpha] & y \\
                w^T  & B[\{k\},\{n+1\}]}
\]
for some $y$ and $w$,
and some $\alpha \subseteq \{1,2, \dots, n\}$ with $|\alpha| = i-1 <r$.
Since the first $|\alpha|=i-1$ columns of $C$ are columns of
the matrix $B[\alpha \cup \{k\}]$, which is nonsingular (because $q_i = \tt A$), the matrix
\[
\mtx{B[\alpha] \\ z^T}
\]
has full rank.
With high probability, $z$ was chosen so that the vector
\[
\mtx{y \\ B[\{k\},\{n+1\}]}
\]
is not in the span of the first $i-1$ columns of $C$;
hence, $C$ is nonsingular, as desired.

\noindent
{\bf Subcase B}: $C$ involves both the last row and the last column of $B'$.

\noindent
Then
\[
C = \mtx{B[\alpha, \beta] & y \\
                w^T  & x^TBx}
\]
for some $\alpha, \beta \subseteq \{1,2, \dots, n\}$ with
$|\alpha| = |\beta| = i-1 <r$ and $|\alpha \cap \beta| = |\alpha|-1$,
where $y = Bx[\alpha]$ and $w = Bx[\beta]$.
Since $B[\alpha, \beta]$ is an $(i-1) \times (i-1)$ submatrix of
an $i \times i$ quasi-principal submatrix of $B$,
and because $q_i = \tt A$,
Lemma \ref{rank when deleting} implies that
$\rank B[\alpha, \beta] \geq i-2$.

\noindent
{\bf Sub-subcase a}: $\rank B[\alpha, \beta] = i-2$.

\noindent
Then with high probability $z$ was chosen so that
$y$ is not in the column space of $B[\alpha, \beta]$ and
$w^T$ is not in the row space of $B[\alpha, \beta]$,
implying that
$\rank C = \rank B[\alpha, \beta]+2 = i$
(see Lemma \ref{rank when appending}), and therefore that
$C$ is nonsingular, as desired.

\noindent
{\bf Sub-subcase b}: $\rank B[\alpha, \beta] = i-1$.

\noindent
In this case, $\rank B[\alpha, \beta]$ has full rank,
meaning that
$y$ is in the column space of $B[\alpha, \beta]$, and that
$w^T$ is in the row space of $B[\alpha, \beta]$.
With high probability, $z$ was chosen so that
the rank of $C$ is larger than that of $B[\alpha, \beta]$, implying that
$\rank C = \rank B[\alpha, \beta]+1 = i$
(see Lemma \ref{rank when appending}), and therefore that
$C$ is nonsingular, as desired.

It follows from Cases 1 and 2 that
$\qpr(B') = q_1q_2 \cdots q_{n-1} \tt SN$, as desired.
\epf

\begin{lem}
\label{AtoSA}
Let $F$ be a field with $\chr F=0$, and let  $B\in\Fnn$ be symmetric.
Assume  that $\qpr(B)=q_1q_2\cdots q_{n-1}{\tt A}$ with $q_k\in\{{\tt A, S}\}$
for $k=1,\dots,n$.
Then there exists a symmetric matrix $B'\in F^{(n+1)\x(n+1)}$ such that
$\qpr(B')=q_1\cdots q_{n-1}{\tt SA}$.
 \end{lem}
\bpf
Suppose $B = [b_1, b_2, \dots, b_n]$ and
let $z \in \SPAN\{b_2, b_3, \dots, b_n\}$ be chosen randomly,
and assume that $z=Bx$.
Let
\[
B' = \mtx{B & z \\
                z^T  & p},
\]
where $p = x^TBx$.
It follows from the proof of Lemma \ref{AtoSN} that
$\qpr(B') = q_1q_2 \cdots q_{n-1} \tt SN$.
%%%%%%%%%%%%%%%%%%%%%%%%%%%%
%%%%%%%%%%%%%%%%%%%%%%%%%%%%
%%%%%%%%%%%%%%%%%%%%%%%%%%%%
%The argument below is the same as the one in Lemma 3.2
%%%%%%%%%%%%%%%%%%%%%%%%%%%%
Define
\[
B'' := \mtx{B & z \\
                z^T  & p + t} =
        \mtx{B & Bx \\
                x^T B & x^TBx + t},
\]
and suppose that $\qpr(B'') = q''_1q''_2 \cdots q''_nq''_{n+1}$.
We now show that, with high probability,
$\qpr(B'')=q_1 \cdots q_{n-1}{\tt SA}$.
Note that $\det B'' = t \det B + \det C = t \det B$.
Then, as $\det B \neq 0$, with high probability, $t$ was chosen so that
$\det B'' \neq 0$, implying that $q''_{n+1} = \tt A$.

Note that
\[
\det B''[\{1, \dots, n\}, \{2,  \dots, n+1\}] =
\det B'[\{1, \dots, n\}, \{2,  \dots, n+1\}] = 0
\]
and
\[
\det B''[\{1, \dots, n\}] = \det B[\{1, \dots, n\}] \neq 0,
\]
implying that
$q''_{n} = \tt S$.

Observe that if $q_j = \tt S$ for some $j$ with $1 \leq j \leq n-1$,
meaning that
$B$ contains both a zero and nonzero
quasi-principal minor of order $j$, then
$B''$ also does; hence,
if  $q_k = \tt S$ for some $k$ with $1 \leq k \leq n-1$, then
$q''_k = \tt S$.

It remains to show that
the following statement holds with high probability:
If $q_j = \tt A$ for some $j$ with $1 \leq j \leq n-1$,
then $q''_j = \tt A$.
Thus, suppose that $q_j = \tt A$ for some $j$ with $1 \leq j \leq n-1$.
Let $B''[\alpha, \beta]$ be a $j \times j$ quasi-principal submatrix of $B''$,
meaning that $|\alpha| = |\beta| = j$ and
$|\alpha| - 1  \leq  |\alpha \cap \beta| \leq |\alpha|$.
Since $q_j = \tt A$,
$\det C[\alpha, \beta] \neq 0$.
If $n+1 \notin \alpha \cup \beta$, then
$\det B''[\alpha, \beta]  = \det C[\alpha, \beta] \neq 0$, as desired.
Now suppose that $n+1 \in \alpha \cup \beta$.
Since $B''$ is symmetric, it suffices to examine the following two cases.

\noindent
{\bf Case i}: $n+1 \in \alpha \cap \beta$.

\noindent
Then
$\det B''[\alpha, \beta] =
t \det B''[\alpha \setminus \{n+1\}, \beta  \setminus \{n+1\}] +
 \det C[\alpha, \beta]$.
Since $\det C[\alpha, \beta] \neq 0$, with high probability,
$t$ was chosen so that $\det B''[\alpha, \beta] \neq 0$.

\noindent
{\bf Case ii}: $n+1 \notin \alpha$ and $n+1 \in \beta$.

\noindent
Suppose that the index of the last row of $B''[\alpha, \beta]$ is $k$,
and that $l \in \alpha $ but $l \notin \beta$.
It is easy to see that
\[
B''[\alpha, \beta] = \mtx{B[\alpha \setminus \{k\}, \beta \setminus \{n+1\}] & y \\
                                            w^T  & B[\{k\},\{n+1\}]}
\]
for some vectors $w, y \in F^{j-1}$.

Since the first $|\alpha|-1=j-1$ columns of $B''[\alpha, \beta] $ are columns of the quasi-principal submatrix
$B[\alpha, \beta \cup \{l\} \setminus \{n+1\}]$,
which is nonsingular (because $q_j = \tt A$),
the matrix
\[
\mtx{B[\alpha \setminus \{k\}, \beta \setminus \{n+1\} ] \\ w^T}
\]
has full rank.
With high probability, $z$ was chosen so that the vector
\[
\mtx{y \\ B[\{k\},\{n+1\}]}
\]
is not in the span of the first $j-1$ columns of $B''[\alpha, \beta] $;
hence, $B''[\alpha, \beta] $ is nonsingular, as desired.

It follows from Cases i and ii that
$\qpr(B'')=q_1 \cdots q_{n-1}{\tt SA}$,
which concludes the proof.
\epf

Given a sequence $t_{i_1}t_{i_2} \cdots t_{i_k}$,
$\overline{t_{i_1}t_{i_2} \cdots t_{i_k}}$ indicates that
the sequence may be repeated as many times as desired
(or it may be omitted entirely).

\begin{thm}\label{allASA}
Let $F$ be a field with $\chr F = 0$.
Any qpr-sequence that does not contain {\tt N} and ends in {\tt A} is
attainable by a symmetric matrix over $F$.
 \end{thm}
\bpf
We proceed by induction.
The sequences {\tt A}, {\tt AA}, and {\tt SA} are attainable
(that is easily checked).
Now assume that all qpr-sequences of length at most $n$ not containing an $\tt N$ and ending with {\tt A} are attainable.
We now show that the sequence $q_1 q_2 \cdots q_n{\tt A}$
from $\{\tt A, S \}$ is attainable.
Since the sequence ${\tt SS\overline{S}A}$ is attainable by the
the direct sum of an identity matrix and the matrix
\[
\mtx{1 & 1 \\
        1 & 0}
\]
(that is easily checked),
we may assume that there exists $j$ with $j \leq n$ such that
$q_j={\tt A}$.
Let $k$ be the largest index for which $q_k={\tt A}$;
hence, $q_{k+1} \cdots q_n=\overline{\tt S}$.
By the induction hypothesis, there is a symmetric matrix $B$ such that
$\qpr(B)=q_1 q_2 \cdots q_k$.
If $k=n$, then $q_1q_2 \cdots q_n{\tt A}$ is attainable by
Lemma  \ref{copyA}.
Thus, assume that $k<n$.
Finally, by applying
Lemma \ref{copyA} followed by
Lemma \ref{AtoSA} $n-k$ times to $q_1 q_2 \cdots q_k$,
the attainability of $q_1 \cdots q_{k}{\tt S\overline{S}A}=
q_1 \cdots q_n{\tt A}$ follows.
\epf

\begin{thm}
\label{allASAN}
Let $F$ be a field with $\chr F = 0$.
Any qpr-sequence  $q_1 q_2 \cdots q_{n}{\tt N}\overline{\tt N}$ with
$q_j \in \{{\tt A, S}\}$ for $j=1,\dots,n$ and $t \geq 1$ copies of $\tt N$
is attainable by a symmetric matrix over $F$.
\end{thm}
\bpf
%The sequence ${\tt S\,\overline{S}\,\overline{N}\,N}$ is attainable (Corollary \ref{SNSN}), so assume there exists $i\le n$ such that $\ell_i={\tt A}$ and let $k$ be the largest index such that $\ell_k={\tt A}$.  It  follows from Theorem \ref{allASA} and Proposition   \ref{copy3}that  the   epr-sequence  $\ell_1\cdots\ell_{n-1}{\tt A}{\tt N}\cdots{\tt N}$  with   $\ell_i\in\{{\tt A,S}\}$  for all $i=1,\dots,n-1$ is attainable. So assume $k\le n-1$ and $\ell_{k+1}=\dots=\ell_n={\tt S}$.
Let $q_1 q_2 \cdots q_n$  be a sequence from $\{\tt A, S \}$ and
$t \geq 1$ be an integer.
First, suppose $q_n={\tt A}$.
By Theorem \ref{allASA}, $q_1 q_2 \cdots q_n$ is attainable.
By applying  Lemma \ref{copyN}  to
$q_1q_2 \cdots q_n$ $t$ times,
the attainability of
$q_1q_2 \cdots q_n{\tt N}\overline{\tt  N}$ follows.

Finally, suppose $q_n={\tt S}$.
Since $q_1q_2 \cdots q_{n-1}{\tt A}$ is
attainable by Theorem \ref{allASA},
applying Lemma \ref{AtoSN} to this sequence shows that
$q_1q_2 \cdots q_{n-1}{\tt SN} =
q_1q_2 \cdots q_n{\tt N} $ is attainable.
Then, by applying Lemma \ref{copyN} $t-1$ times to
$q_1q_2 \cdots q_n{\tt N}$, the attainability of
$q_1q_2 \cdots q_n{\tt N}\overline{\tt  N}$ follows.
\epf

Corollary \ref{cor: necessary condition} and
Theorems \ref{allASA} and \ref{allASAN} can now be used to obtain a complete characterization of all the qpr-sequences that are attainable by symmetric matrices over a field of characteristic $0$:

\begin{thm}
\label{thm: main result}
The qpr-sequence $q_1q_2 \cdots q_n$ is attainable by a
symmetric matrix over a field of characteristic 0
if and only if
the following hold:
\ben
\item[(i)] $q_n \neq \tt S$.
\item[(ii)] Neither $\tt NA$ nor $\tt NS$ is a subsequence of $q_1q_2 \cdots q_n$.
\een
\end{thm}

Theorem \ref{thm: main result} does not
hold in general for fields of nonzero characteristic:
Consider a qpr-sequence of the form
$\tt AA \overline{A}N\overline{N}$.
Any matrix with this qpr-sequence must, obviously,
have epr-sequence $\tt AA \overline{A}N\overline{N}$
(see Observation \ref{obs: qpr with A and N}),
which, by \cite[Theorem 3.8]{XMR-Char 2}, is not attainable by
symmetric matrices over the prime field of order 2
(i.e., the integers modulo 2), implying that the qpr-sequence
$\tt AA \overline{A}N\overline{N}$ is not attainable by
symmetric matrices over this field (which is of characteristic 2).

%%%%%%%%%%%%%%%%%%%%%%%%%%%%%%%%%%%%%%%
%%%%%%%%%%%%%%%%%%%%%%%%%%%%%%%%%%%%%%%
%%%%%%%%%%%%%%%%%%%%%%%%%%%%%%%%%%%%%%%
%%%%%%%%%%%%%%%%%%%%%%%%%%%%%%%%%%%%%%%
%%%%%%%%%%%%%%%%%%%%%%%%%%%%%%%%%%%%%%%
%%%%%%%%%%%%%%%%%%%%%%%%%%%%%%%%%%%%%%%
%%%%%%%%%%%%%%%%%%%%%%%%%%%%%%%%%%%%%%%
%%%%%%%%%%%%%%%%%%%%%%%%%%%%%%%%%%%%%%%
\subsection*{Acknowledgments}
$\null$
\indent
The second author expresses his gratitude to the Department of Mathematics at Iowa State University, for its wonderful hospitality while he was a graduate student there, when the research presented here was conducted.

%%%%%%%%%%%%%%%%%%%%%%%%%%%%%%%%%%%%%%%
%%%%%%%%%%%%%%%%%%%%%%%%%%%%%%%%%%%%%%%
%%%%%%%%%%%%%%%%%%%%%%%%%%%%%%%%%%%%%%%
%%%%%%%%%%%%%%%%%%%%%%%%%%%%%%%%%%%%%%%
%%%%%%%%%%%%%%%%%%%%%%%%%%%%%%%%%%%%%%%
%%%%%%%%%%%%%%%%%%%%%%%%%%%%%%%%%%%%%%%
%%%%%%%%%%%%%%%%%%%%%%%%%%%%%%%%%%%%%%%
%%%%%%%%%%%%%%%%%%%%%%%%%%%%%%%%%%%%%%%

\end{document}